\documentclass{amsart}

\usepackage{amssymb}
\usepackage[colorlinks=false, urlcolor=blue, final]{hyperref}
\usepackage[all]{xy}
\usepackage{enumitem}

\newtheorem{thm}{Theorem}[section]
\newtheorem{lem}[thm]{Lemma}
\newtheorem{cor}[thm]{Corollary}
\newtheorem{prop}[thm]{Proposition}

\theoremstyle{definition}

\newtheorem{ex}[thm]{Example}
\newtheorem{rmk}[thm]{Remark}
\allowdisplaybreaks

\numberwithin{equation}{section}

\theoremstyle{remark}

\newtheorem{nota}[thm]{Notation}

\DeclareMathOperator{\lt}{{lt}}
\DeclareMathOperator{\id}{id}
\DeclareMathOperator{\supp}{supp}
\DeclareMathOperator{\spn}{span}
\DeclareMathOperator{\Rep}{Rep}

\newcommand{\T}{\mathbb{T}}

\newcommand{\C}{\mathbb{C}}
\newcommand{\Z}{\mathbb{Z}}
\newcommand{\units}{G^{(0)}}
\newcommand{\inv}{^{-1}}
\newcommand{\linner}[3]{{{}_{{}_{#1}}}\!\left\langle #2, #3\right\rangle}
\newcommand{\rinner}[3]{\left\langle #1, #2\right\rangle_{{}_{#3}}}

\begin{document}

\title[The $C^*$-algebras of groupoid extensions]
{\boldmath{Decomposing the $C^*$-algebras of groupoid extensions}}

\author[J. H. Brown]{Jonathan H. Brown}
\author[A. an Huef]{Astrid an Huef}
\address{Department of Mathematics and Statistics\\
University of Otago\\
Dunedin 9054\\
New Zealand}
\email{jbrown@maths.otago.ac.nz}
\email{astrid@maths.otago.ac.nz}

\date{April 28, 2011, with minor revisions October 8, 2011 and May 3, 2012.}

\begin{abstract} We decompose the full and reduced $C^*$-algebras of an extension of a  groupoid by the circle into a direct sum of twisted groupoid $C^*$-algebras. 
\end{abstract}

\thanks{We thank Iain Raeburn and Dana Williams for helpful discussions}
\commby{Marius Junge}
\maketitle

\section{Introduction}  Let $H$ be a compact group and denote by $\hat H$ the collection of equivalence classes of the irreducible unitary representations of $H$.  The Peter-Weyl Theorem implies first, that every irreducible unitary representation of $H$ is finite-dimensional, and, second, that the left-regular representation $\lambda$ of $H$ on $L^2(H)$ is unitarily equivalent to the direct sum $
\oplus_{[U]\in\hat H}d_U\cdot  U$, where $d_U$ is the dimension of $U$. The $C^*$-algebra $C^*(H)$ of $H$ is universal for the unitary representations of $H$, which means roughly that the unitary representations $U$ of $H$ are in one-to-one correspondence with the nondegenerate representations $\pi_U$ of $C^*(H)$ on the Hilbert space  of $U$. Since $H$ is compact, the left-regular representation $\pi_\lambda$ is an isomorphism, and the reduced $C^*$-algebra $C^*_r(H):=\pi_\lambda(C^*(H))$ coincides with $C^*(H)$.  So in the $C^*$-setting, the Peter-Weyl theorem says that $C^*(H)=C^*_r(H)$ is a direct sum $\oplus_{[U]\in\hat H} M_{d(U)}(\C)$ of matrix algebras.

A similar result holds for extensions of locally compact groups $H$ by the circle $\T$.  Let $\omega:H\times H\to \T$ be a continuous $2$-cocycle. We associate  two $C^*$-algebras to the pair $(H, \omega)$.  The first is  the twisted group $C^*$-algebra $C^*(H,\omega)$ which is universal for the $\omega$-representations of $H$. For the second, equip $H^\omega:=\T\times H$  with the product topology and  multiplication $(s,\eta)(t,\gamma)=(st\omega(\eta,\gamma),\eta\gamma)$; then $H^\omega$ is a locally compact group and has a $C^*$-algebra. It follows from \cite{DPW89}, for example, that $C^*(H^\omega)$ is isomorphic to the direct sum $\oplus_{n\in\Z} C^*(H, \omega^n)$ of twisted group $C^*$-algebras (see also \cite[Corollary~3]{Pah93}).
In this paper we generalize this latter result to locally compact Hausdorff groupoids.  

Let $G$ be a locally compact Hausdorff groupoid and $\omega:G^{(2)}\to\T$ a continuous $2$-cocycle on the set $G^{(2)}$ of composable pairs in $G$.  We show that the $C^*$-algebra of the extension  $G^\omega$ is isomorphic to the  direct sum $\oplus_{n\in\Z} C^*(G, \omega^n)$ of twisted groupoid  $C^*$-algebras, and that this isomorphism factors through to the reduced $C^*$-algebras (see Theorems~\ref{thm-dec} and \ref{thm-red dec}).  The full twisted groupoid $C^*$-algebras have been used in \cite{MW92}, \cite{MRW96} and \cite{CH10} to characterize when groupoid $C^*$-algebras have continuous trace or bounded trace, and their non-selfadjoint subalgebras have been studied in  \cite{MS89}. The reduced twisted groupoid $C^*$-algebras  appear  as the $C^*$-algebras with diagonal subalgebras in \cite{Kum86} and as the $C^*$-algebras with Cartan subalgebras in \cite{Ren08}.  For example,  if $A$ is a $C^*$-algebra with diagonal subalgebra $B$ then Kumjian's Theorem~3.1 of \cite{Kum86} implies that  there exists a principal \'etale groupoid $G$ and an extension of $G$ by $\T$  implemented by a (possibly Borel) cocycle $\omega$  such that $A$ is isomorphic to $C^*_r(G, \overline{\omega})$, and the isomorphism maps the diagonal $B$ to a diagonal in $C^*_r(G, \overline{\omega})$.  A similar result holds for Cartan subalgebras, except there the groupoid $G$ may be only topologically principal \cite[\S5]{Ren08}.

The main theorems of this paper provide a general framework for investigating twisted groupoid $C^*$-algebras using the literature on  the non-twisted case.  For example, suppose that  $G$ is principal. Then we  deduce in Proposition~\ref{prop-similar}  that $C^*(G)$ has continuous trace if and only if $C^*(G^\omega)$ has continuous trace, and if $C^*(G)$ has continuous trace then so does $C^*(G,\omega)$.  See Proposition~\ref{prop-similaragain} for more results along these lines.  We also deduce in  Corollary~\ref{cor amen}  that if a groupoid $G$ is amenable then $C^*(G,\omega)$ and $C^*_r(G,\omega)$ are isomorphic.

\section{Preliminaries}
Throughout, $G$ is a  second-countable, locally compact, Hausdorff groupoid with Haar sytstem $\{\lambda^{u}\}_{u\in\units}$.    We denote by $\lambda_u$ the image of $\lambda^u$ under inversion.  We write
 $G^{(0)}$ for the unit space of $G$,
$r=r_G,s=s_G:G\to G^{(0)}$  for the range and source maps $r_G(\gamma)=\gamma\gamma^{-1}$ and
$s_G(\gamma)=\gamma^{-1}\gamma$, respectively, and $G^{(2)}:=\{(\gamma,\eta):s_G(\gamma)=r_G(\eta)\}$ for the set of composable pairs.  
  
\subsection{Groupoid extensions} Let $\omega:G^{(2)}\to\T$  be a continuous   $2$-cocycle, so  that $\omega$ satisfies the cocycle  identity \[\omega(\gamma,\eta)\omega(\gamma\eta,\xi)=\omega(\eta,\xi)\omega(\gamma,\eta\xi).\]  We will assume throughout  that $\omega$ is normalized in the sense that  
\[
\omega(r_G(\gamma),\gamma)=1=\omega(\gamma,s_G(\gamma))\text{\ and \ }\gamma\in G;
\]  since every $2$-cocycle is cohomologous to a normalized one\footnote{The coboundary implementing the equivalence is the image under the boundary map of the function $b(\gamma)=\omega(r_G(\gamma),\gamma)$.}  and because  the associated $C^*$-algebras  depend only on the class  of the $2$-cocycle (see \cite[Proposition~II.1.2]{Ren80}), there is no loss of generality. Following~\cite[page~73]{Ren80} we denote by $G^\omega$ the extension of $G$ by $\T$ defined by $\omega$:  thus   $G^\omega$ is the groupoid $\T\times G$  with the product topology, with  range and source  maps $r_{G^\omega}(t,\gamma)=(1,r_G(\gamma))$ and $s_{G^\omega}(t,\gamma)=(1,s_G(\gamma))$, multiplication $(s,\eta)(t,\gamma)=(st\omega(\eta,\gamma),\eta\gamma)$ and inverse $(t,\gamma)^{-1}=(t^{-1}\omega(\gamma,\gamma^{-1})\inv,\gamma^{-1})$.  We  identify the unit space of $G^\omega$ with $G^{(0)}$ via $(1,u)\mapsto u$.  We say that $G^\omega$ is the \emph{groupoid extension  associated to $(G,\omega)$}.  When we want to emphasize the product nature of $G^\omega$ we will denote it by $\T\times_\omega G$.

In order to reconcile our work with the literature, suppose that 
\begin{equation}\label{eq ext}\units\to \T\times\units\stackrel{i}{\hookrightarrow} E\stackrel{j}{\to} G:=E/\T\to \units
\end{equation} 
is an extention of topological groupoids such that $i$ induces a free action of $\T$ on $E$ by $t\cdot\gamma=i(t, r_E(\gamma))\gamma$ for $\gamma\in E$ and $t\in\T$.
In \cite[page~131]{MW92}, Muhly and Williams discuss a correspondence between extensions $E$ and Borel $2$-cocycles defined using a Borel cross section of $j$.  They show that the $2$-cocycle $\omega$ associated to an extension $E$ 
is continuous if and only if there exists a continuous section of $j$, and then $E$ is topologically isomorphic to $G^\omega$. (If the cocycle $\omega$ is not continuous, then $G^\omega$ is a Borel groupoid and the topological groupoid $E$ is Borel isomorphic to $G^\omega$.) When $E=G^\omega$ for a continuous $\omega$,  the inclusion $i$ in \eqref{eq ext} is  $i(t,u)=(t,u)$ and $j$ is the projection onto $G$.  Furthermore, since $\omega$ is normalized, the action of $\T$ on $G^\omega$ induced by $i$ is  $s\cdot(t,\gamma)=(st,\gamma)$. Our main reason for restricting our attention to extensions associated to continuous cocycles is that we use the continuity in an apparently essential way in Lemma~\ref{lem: density of poly}.

\begin{ex}\label{ex lots of co}
It is quite easy to construct groupoids $G$ with non-trivial continuous $2$-cocycles $\omega$, and hence there are  many non-trivial extensions $G^\omega$ as described above.  For example, take $X=S^3$ and recall that  the \v{C}ech cohomology group $H^3(X,\Z)$ is non-trivial. Since $H^3(X,\Z)$ is isomorphic to the second sheaf cohomology group $H^2(X,\mathcal{S})$ (see, for example, \cite[Theorem~4.42]{tfb}), there exists a non-trivial cocycle  $\lambda:=\{\lambda_{ijk}:U_{ijk}\to\T\}$ where the $U_{ijk}$ are the  triple overlaps of  an open cover $\{U_i\}$ of $X$.  Define \[\Psi:\bigcup_i U_i\times\{i\}\to X\text{\ by $\Psi(x,i)=x$ for $x\in U_i$}.\]
Then $R(\Psi):=\{ ((x,i),(x,j)):x\in U_i\cap U_j \}$ becomes a groupoid with range and source maps given by  $r_{R(\Psi)}((x,i),(x,j))=(x,i)$ and  $s_{R(\Psi)}((x,i),(x,j))=(x,j)$, multiplication defined by $((x,i),(x,j))((x,j),(x,k))=((x,i),(x,k))$ and inverse $((x,i),(x,j))^{-1}=((x,j),(x,i))$. Now define $\omega_\lambda:R(\Psi)^{(2)}\to\T$ by 
\[
\omega_\lambda \big(((x,i),(x,j)),((x,j),(x,k))\big)=\lambda_{ijk}(x)
\]
for $x\in U_{ijk}$.
It is straightforward to check that $\omega_\lambda$ is a non-trivial continuous $2$-cocyle (the  $\lambda_{ijk}:U_{ijk}\to \T$ are continuous by definition, and $\omega_\lambda$ is a coboundary if and only if $\lambda$ is).
\end{ex}

Recall that a groupoid $G$ is {\emph{principal}} if the map $\Phi:\gamma\mapsto (r_G(\gamma),s_G(\gamma))$ is injective  and is \emph{proper} if  $\Phi$ is proper.  We say $G$ is  \emph{transitive} if given $u,v\in \units$ there exists $\gamma\in G$ such that $r_G(\gamma)=u$ and $s_G(\gamma)=v$.

\begin{rmk} Although Example~\ref{ex lots of co} shows that there are many examples of nontrivial continuous $2$-cocycles, principal transitive groupoids have only trivial ones.  To see this, let $G$ be a principal and transitive groupoid,  and let $\omega$  be a normalized $2$-cocycle on $G$.  Pick $u\in \units$, and let $b:\gamma\mapsto \omega(\gamma,\alpha_\gamma)$ where $\alpha_\gamma$ is the unique element such that $r_G(\alpha_\gamma)=s_G(\gamma)$ and $s_G(\alpha_\gamma)=u$. Then $\omega(\gamma,\eta)=b(\gamma)b(\eta)\overline{b(\gamma\eta)}$ and thus $\omega$ is a coboundary.
\end{rmk}


\subsection{The $C^*$-algebras}\label{subsec-c*alg}
Let  $E$ be a second-countable, locally compact groupoid with left Haar system $\beta$, and $\sigma:E^{(2)}\to \T$ a continuous, normalized $2$-cocycle.  For $f,g\in C_c(E)$, the formulas
\begin{align*}
f*g(\gamma)&=\int_E  f(\eta)g(\eta^{-1}\gamma)\sigma(\eta,\eta^{-1}\gamma)\,\, d\beta^{r_E(\gamma)}(\eta)\text{\ \ and \ \ }
f^*(\gamma)=\overline{f(\gamma^{-1})}\,\overline{\sigma(\gamma,\gamma^{-1})}
\end{align*}
define a convolution and involution on $C_c(E)$.
These operations  make $C_c(E)$ into a $*$-algebra, denoted by $C_c(E,\sigma)$.  We denote by $\Rep(C_c(E,\sigma))$ the set of  Hilbert-space representations $\rho:C_c(E,\sigma)\to B(\mathcal{H})$  that are continuous for the inductive limit topology on $C_c(E)$ and the weak operator topology on $B(\mathcal{H})$.   Then
\begin{equation}\label{eq-norm}\|f\|=\sup\{\|\rho(f)\|:\rho\in \Rep(C_c(E,\sigma))\}\end{equation}
is finite and defines a pre-$C^*$-norm on $C_c(E,\sigma)$; the \emph{twisted groupoid $C^*$-algebra} $C^*(E,\sigma)$ is  defined as the  completion of $C_c(E,\sigma)$ in this norm.  (All of this is non-trivial.  If $\rho\in\Rep(C_c(E,\sigma))$ then $\rho$ is the integrated form of a unitary representation of $G$ by  \cite[Th\`{e}or\'{e}me~4.1(i)]{Ren87}, and then    $\rho$ is bounded in Renault's $I$-norm \cite[Proposition~II.1.7]{Ren80}.  Representations bounded by the $I$-norm are continuous in the inductive limit topology. It now  follows that \eqref{eq-norm} defines a norm by \cite[Proposition~II.1.11]{Ren80} and that the definition of $C^*(E,\sigma)$ above coincides with the one in \cite[Definition~II.1.12]{Ren80}.)  If the cocycle is identically $1$ then we write $C^*(E)$ for $C^*(E, 1)$ and call it the \emph{groupoid $C^*$-algebra} of $E$.  

Let  $\tau$ be normalized left Haar measure on $\T$; we will denote $d\tau(t)$ by $dt$. 
Let $\omega$ be a continuous $2$-cocycle on $G$ and let $G^\omega$ be the associated groupoid extension.  Since $G^\omega=\T\times_\omega G$ has the product topology,  the product measures $\{\tau\times \lambda^u:u\in G^{(0)}\}$  define a  Haar system on the extension $G^\omega$.   
For fixed $n\in\Z$, let \[C_c(G^\omega,n)=\{f\in C_c(G^\omega):f(s\cdot(t,\gamma))=s^{-n}f(t,\gamma)\}.\] As above, we denote by $\Rep(C_c(G^\omega,n))$ the set of Hilbert-space representations   $\rho:C_c(G^\omega,n)\to B(\mathcal{H})$ that are continuous for the inductive limit topology on $C_c(G^\omega,n)$  and the weak operator topology on $B(\mathcal{H})$.  Then  $C_c(G^\omega,n)$ is a $*$-subalgebra of $C_c(G^\omega)$, and, as in \cite[\S5]{Ren87} and \cite[page~130]{MW92}, the $C^*$-algebra $C^*(G^\omega,n)$ is the completion of $C_c(G^\omega,n)$ in the norm $\|f\|=\sup\{\|\rho(f)\|:\rho\in \Rep C_c(G^\omega,n)\}$. (Again, this is non-trivial: Corollaire~4.8 of \cite{Ren87} implies that this indeed defines a norm bounded by the $I$-norm.) The $C^*$-algebra $C^*(G^\omega,n)$ was studied in \cite[\S1]{Ren87},  and, when $n=-1$, in \cite{MW92}.

\begin{rmk} Let $f,g\in C_c(G^\omega,n)$.   Since the Haar system $\{\tau\times\lambda^u\}$ on $G^\omega$ is pulled back from the one on $G$ and $\tau$ is normalized, the convolution $f*g$ can be written as an integral over $G$: a direct calculation shows that for any $s\in\T$,
\begin{align*} f*g(t,\gamma)&:=\int_G\int_\T f(r,\eta)g((r,\eta)\inv(t,\gamma))\, dr \, d\lambda^{r_G(\gamma)}(\eta)\\
&=\int_G\int_\T r^{-n}r^n f(1,\gamma)g( t\omega(\eta,\eta\inv)\inv\omega(\eta\inv,\gamma),\eta\inv\gamma)\, dr \, d\lambda^{r_G(\gamma)}(\eta)\\
&=\int_G  f(s,\eta)g((s,\eta)\inv (t,\gamma))\, \, d\lambda^{r_G(\gamma)}(\eta).
\end{align*}\end{rmk}


\section{Decomposing the $C^*$-algebra of a groupoid extension}\label{sec-dec}
Throughout $\omega:G^{(2)}\to\T$ is a continuous normalized $2$-cocycle, and  $G^\omega$ is the groupoid extension associated to $(G,\omega)$.  Note that $\omega^n$ is also a continuous $2$-cocycle.  
The goal of this section is to prove that $C^*(G^\omega)$ is isomorphic to a direct sum of twisted groupoid $C^*$-algebras $C^*(G, \omega^n)$.   We start by proving  that $C^*(G^\omega,n)$ is a quotient of $C^*(G^\omega)$ and is isomorphic to $C^*(G,\omega^n)$.

\begin{lem}\label{lem-algberas} Let $G^\omega$ be the groupoid extension associated to $(G,\omega)$. Fix $n\in\Z$.
\begin{enumerate}
\item\label{item: chi homo}\cite[Lemma~3.3]{Ren87}
Define $\chi_n:C_c(G^\omega)\rightarrow C_c(G^\omega,n)$ by
$$\chi_n(f)(t,\gamma):=\int_\T f(s\cdot (t,\gamma))s^{n} \, ds=\int_\T f(st,\gamma)s^{n} \, ds.$$
Then $\chi_n$ is a $*$-homomorphism continuous with respect to  the inductive limit topologies, and  extends to a $*$-homomorphism $\chi_n:C^*(G^\omega)\to C^*(G^\omega,n)$ such that  $\chi_n(f)=f$ for $f\in C_c(G^\omega,n)$. In particular, $\chi_n$ is a quotient map.
\item \label{item:isotoRenault}
Let $\phi_n: C_c(G^\omega,n)\to C_c(G, \omega^n)$ be the map $\phi_n(f)(\gamma)=f(1,\gamma)$ for $\gamma\in G$. Then $\phi_n$ extends to a $*$-isomorphism of $C^*(G^\omega,n)$ onto $C^*(G,\omega^n)$.
\end{enumerate}
\end{lem}

\begin{proof} Part~\eqref{item: chi homo} is  \cite[Lemma~3.3]{Ren87} (see also \cite[Proposition II.1.22]{Ren80} for a detailed proof of  the case $n=1$).

(\ref{item:isotoRenault}) It suffices to show that $\phi_n: C_c(G^\omega,n)\to C_c(G, \omega^n)$ is a continuous bijective $*$-homomorphism with a continuous inverse. For then $\phi_n$ and $\phi_n\inv$ extend to $*$-homomorphisms $\phi_n :C^*(G^\omega,n)\to C^*(G,\omega^n)$ and  $\phi_n\inv:C^*(G,\omega^n)\to C^*(G^\omega,n)$, and by continuity $\phi_n\circ\phi_n\inv=\id$ and $\phi_n\inv\circ\phi_n=\id$, giving that  $\phi_n$ is an isomorphism.

To see that $\phi_n$ is a homomorphism, we first need a calculation with cocycles. Let $\eta,\gamma\in G$ with $r_G(\gamma)=r_G(\eta)$. Since $\omega$ is  normalized, we have 
\begin{gather*}
1=\omega(r_G(\eta^{-1}\gamma),\eta^{-1}\gamma)=\omega(s_G(\eta),\eta^{-1}\gamma)=\omega(\eta^{-1}\eta,\eta^{-1}\gamma) \quad{\text{and}}\\ 
\omega(\eta,\eta\inv)=\omega(\eta,\eta\inv)\omega(\eta\eta\inv,\eta)=\omega(\eta\inv,\eta)\omega(\eta,\eta\inv\eta)=\omega(\eta\inv,\eta).
\end{gather*}
Thus  
\[\omega(\eta\inv,\eta)=\omega(\eta\inv,\eta)\omega(\eta\inv\eta,\eta\inv\gamma)=\omega(\eta,\eta\inv\gamma)\omega(\eta\inv,\gamma),
\]
and 
it follows that  \begin{equation}\label{eqcocycle}\overline{\omega(\eta,\eta\inv)}\omega(\eta\inv,\gamma)=\overline{\omega(\eta,\eta\inv\gamma)}.\end{equation}
So, for $f,g\in C_c(G^\omega,n)$, 
\begin{align*} \phi_n(f*g)(\gamma)&=f*g(1,\gamma)=\int_G f(1,\eta)g((1,\eta)\inv(1,\gamma))\, d\lambda^{r_G(\gamma)}(\eta)\\
&=\int_G f(1,\eta)g(\overline{\omega(\eta,\eta\inv)}\omega(\eta\inv,\gamma),\eta\inv\gamma)\, d\lambda^{r_G(\gamma)}(\eta)\\
&=\int_G f(1,\eta)g(\overline{\omega(\eta,\eta\inv\gamma)},\eta\inv\gamma)\, d\lambda^{r_G(\gamma)}(\eta)\quad\text{(using \eqref{eqcocycle})}\\
&=\int_G f(1,\eta)g(1,\eta\inv\gamma)\omega(\eta,\eta\inv\gamma)^n \, d\lambda^{r_G(\gamma)}(\eta)\\
&=\phi_n(f)*\phi_n(g)
\end{align*}
and 
\begin{align*}\phi_n(f^*)(\gamma)
&=\overline{f((1,\gamma)\inv)}
=\overline{f(\omega(\gamma,\gamma\inv)\inv,\gamma\inv)}\\
&=\overline{f(1,\gamma\inv)}\overline{\omega(\gamma,\gamma\inv)^n}
=\phi_n(f)^*(\gamma).\end{align*}
So  $\phi_n$ is a $*$-homomorphism.
To see $\phi_n$  is injective on $C_c(G^\omega, n)$, suppose $f(1,\gamma)=g(1,\gamma)$ for all $\gamma\in G$. Then for all $t\in \T$, 
$$f(t,\gamma)=t^{-n}f(1,\gamma)=t^{-n}g(1,\gamma)=g(t,\gamma)$$ 
and thus $f=g$.  To see $\phi_n$ is onto $C_c(G,\omega^n)$, let $f\in C_c(G)$ and note that $(t,\gamma)\mapsto t^{-n}f(\gamma)$ is in $C_c(G^\omega,n)$, and  $\phi_n$ sends it back to $f$. So $\phi_n :C_c(G^\omega,n)\to C_c(G, \omega^n)$ is a bijection.

If $F_i\to F$ in the inductive limit topology on $C_c(G^\omega)$, then $F_i(1,\cdot)\to F(1,\cdot)$ uniformly on a fixed compact set as well. Thus $\phi_n$ is continuous for the inductive limit topology on $C_c(G^\omega,n)$ and extends to a $*$-homomorphism of the $C^*$-algebras.  Similarly, if $f_i\rightarrow f\in C_c(G)$ in the inductive limit topology, then  $|t^{-n}f_i(\gamma)-t^{-n}f(\gamma)|\leq|f_i(\gamma)-f(\gamma)|$ is eventually small, so that $\phi_n\inv(f_i)\rightarrow \phi_n\inv(f)$ in the  inductive limit topology as well. As outlined at the beginning of the proof, this implies that $\phi_n$ extends to an isomorphism of $C^*(G^\omega,n)$ onto $C^*(G,\omega^n)$.
\end{proof}

Define $\Upsilon_n:=\phi_{n}\circ \chi_{n}: C^*(G^\omega)\to C^*(G,\omega^n)$; then
\[
\Upsilon_n(F)(\gamma)=\int_\T F(t, \gamma)t^{n}\, \, dt\quad\text{for $F\in C_c(G^\omega)$}.
\]

\begin{thm}\label{thm-dec} Let $G$ be a second-countable,  locally compact Hausdorff groupoid with a Haar system $\lambda$.  Let $\omega:G^{(2)}\to\T$ be a continuous  $2$-cocycle and let $G^\omega$ be the groupoid extension associated to $(G,\omega)$.  Then  the map $\Upsilon:C_c(G^\omega)\to \oplus_{n\in\Z} C_c(G, \omega^n)$ defined by $F\mapsto (\Upsilon_n(F))$ extends to an isomorphism of $C^*(G^\omega)$ onto $\oplus_{n\in\Z} C^*(G, \omega^n)$. 
\end{thm}

To prove Theorem~\ref{thm-dec} we first prove that the subalgebra $I_n:=\overline{C_c(G^\omega,n)}^{\|\cdot\|_{C^*(G^\omega)}}$ is an ideal of $C^*(G^\omega)$ which is isomorphic to $C^*(G^\omega,n)$, and second,  that $C^*(G^\omega)$ is the (internal) direct sum of the $I_n$.

\begin{lem}\label{lem:chi iso} Let $G^\omega$ be the groupoid extension associated to $(G,\omega)$. 
\begin{enumerate}
\item\label{item: chi iso} For $f\in C_c(G^\omega,n)$, $\|f\|_{C^*(G^\omega)}=\|f\|_{C^*(G^\omega,n)}$.
\item\label{former cor} The map $\chi_n:C_c(G^\omega,n)\subset C_c(G^\omega)\rightarrow C_c(G^\omega,n)$ extends to an isometry $\chi_n$ of the subalgebra $I_n$ of $C^*(G^\omega)$ onto $C^*(G^\omega,n)$. 
\item\label{item3:chi iso} The quotient map $\chi_n: C^*(G^\omega)\to C^*(G^\omega,n)$ is identically zero on $I_m$ if $n\neq m$.
\end{enumerate}
\end{lem}

\begin{proof} \eqref{item: chi iso} Fix $f\in C_c(G^\omega,n)$. If $\pi\in\Rep(C_c(G^\omega,n))$ then by   Lemma~\ref{lem-algberas}\eqref{item: chi homo}, $\pi\circ \chi_n\in \Rep(C_c(G^\omega))$.  Since $f=\chi_n(f)$ we have 
\begin{align*}
\|f\|_{C_c(G^\omega,n)}&=\sup\{\|\pi(f)\|:\pi\in \Rep(C_c(G^\omega,n))\}\\
&=\sup\{\|\pi\circ\chi_n(f)\|:\pi\in \Rep(C_c(G^\omega,n))\}\\
&\leq \|f\|_{C^*(G^\omega)}.
\end{align*}

Conversely, if $\rho\in\Rep(C_c(G^\omega))$  then $\rho|_{C_c(G^\omega,n)}$ is also continuous in the inductive limit topology on $C_c(G^\omega,n)$.  Fix $\epsilon>0$.  Pick a representation $\rho$ of $C_c(G^\omega)$ such that $\|f\|_{C^*(G^\omega)}<\|\rho(f)\|+\epsilon$.  Then
$$\|f\|_{C^*(G^\omega)}<\|\rho(f)\|+\epsilon=\|\rho|_{C_c(G^\omega,n)}(f)\|+\epsilon\leq\|f\|_{C^*(G^\omega,n)}+\epsilon.$$
Thus $\|f\|_{C^*(G^\omega)}\leq\|f\|_{C^*(G^\omega,n)}$, and  $\|f\|_{C^*(G^\omega)}=\|f\|_{C^*(G^\omega,n)}$ as desired.

\eqref{former cor} Fix $g\in I_n$. Let $\{f_i\}\subset C_c(G^\omega,n)$ be a sequence converging to $g$.
By  \eqref{item: chi iso},  $\|f_i\|_{C^*(G^\omega)}=\|\chi_n(f_i)\|_{C^*(G^\omega,n)}$, and hence $\|g\|_{C^*(G^\omega)}=\|\chi_n(g)\|_{C^*(G^\omega,n)}$.
So $\chi_n$ is isometric on the subalgebra $I_n$ of $C^*(G^\omega)$.  Furthermore, $\chi_n|_{I_n}$ is onto since $C_c(G^\omega,n)$ is dense in $C^*(G^\omega,n)$.  So $\chi_n|_{I_n}$ is an isomorphism.

\eqref{item3:chi iso} This is a direct calculation.
\end{proof}

\begin{lem}\label{lem-Im} Let $G^\omega$ be the groupoid extension associated to $(G,\omega)$.  For each $n\in \Z$,  $I_n$ is an ideal in $C^*(G^\omega)$. Furthermore, $I_n I_m=\{0\}$   if $n\neq m$.\end{lem}

\begin{proof}Let $f\in C_c(G^\omega)$ and $g\in C_c(G^\omega,n)$. Then
\begin{align*} f*g(s\cdot (t,\gamma))&=\int_G\int_\T f(r,\eta)g((r,\eta)\inv(st,\gamma))\, dr \, d\lambda^{r_G(\gamma)}(\eta)\\
&=\int_G\int_\T f(r,\eta)g(sr\inv t \overline{\omega(\eta,\eta\inv)}\omega(\eta\inv,\gamma),\eta\inv\gamma)\, dr \, d\lambda^{r_G(\gamma)}(\eta)\\
&=\int_G\int_\T f(r,\eta)g(s\cdot((r,\eta)\inv(t,\gamma)))\, dr \, d\lambda^{r_G(\gamma)}(\eta)\\
&=s^{-n}\int_G\int_\T f(r,\eta)g((r,\eta)\inv(t,\gamma))\, dr \, d\lambda^{r_G(\gamma)}(\eta)\\
&=s^{-n} f*g(t,\gamma).
\end{align*}
Thus $f*g\in C_c(G^\omega,n)\subset I_n$.  Since $C_c(G^\omega, n)$ is closed under involution $g*f\in I_n$ as well. 
Since $I_n$ is closed the above calculations show that $I_n$ is an ideal in $C^*(G^\omega)$.

To see that $I_m I_n=\{0\}$ unless $n=m$,
 let $f\in C_c(G^\omega,m),g\in C_c(G^\omega,n)$. Then
\begin{align*} f*g(t,\gamma) &=\int_G\int_\T f(r,\eta)g(r\inv t \overline{\omega(\eta,\eta\inv)}\omega(\eta\inv,\gamma),\eta\inv\gamma)\, dr \, d\lambda^{r_G(\gamma)}(\eta)\\
&=\int_G\int_\T r^{-m}f(1,\eta)r^{n}g( t \overline{\omega(\eta,\eta\inv)}\omega(\eta\inv,\gamma),\eta\inv\gamma)\, dr \, d\lambda^{r_G(\gamma)}(\eta)\\ 
&=\int_Gf(1,\eta)g( t \overline{\omega(\eta,\eta\inv)}\omega(\eta\inv,\gamma),\eta\inv\gamma)\, d\lambda^{r_G(\gamma)}(\eta)\int_\T r^{n-m}\, dr. \\
&= \begin{cases} \int_Gf(1,\eta)g( t \overline{\omega(\eta,\eta\inv)}\omega(\eta\inv,\gamma),\eta\inv\gamma)\, d\lambda^{r_G(\gamma)}(\eta) & \text{if $m=n$}\\
0 & \text{otherwise}.\end{cases}\qedhere
\end{align*}
\end{proof}

\begin{nota}For $f\in C_c(G)$, $\psi\in C(\T)$, denote by $\psi\otimes f$ the function $(t,\gamma)\mapsto \psi(t)f(\gamma)$.  In particular, for fixed $m$, we write $s^{m}\otimes f$ for  the function $(t,\gamma)\mapsto t^{m}f(\gamma)$ in $I_{-m}$.\end{nota}

\begin{lem}\label{lem: density of poly} Let $G^\omega$ be the groupoid extension associated to $(G,\omega)$.  The  $\spn\{s^m\otimes f: m\in \Z, f\in C_c(G)\}$ is dense in $C_c(G^\omega)$ in the inductive limit topology.\end{lem}

\begin{proof} Fix $F\in C_c(G^\omega)$ and $\epsilon>0$. Let $U_1$ and $U_2$ be open, relatively compact neighborhoods in $\T$ and $G$, respectively, such that $\supp F\subset U_1\times U_2$.  Because $\omega$ is continuous, $G^\omega=\T\times_\omega G$ has the product topology, and  the map $t\mapsto F(t,\cdot)$ is in $C_c(\T, C_c(G))$. So the support of $t\mapsto F(t,\cdot)$ is contained in $U_1$.  For each $t\in U_1$ let 
\[W_t:=\{s\in \T: \|F(s,\cdot)-F(t,\cdot)\|_\infty<\epsilon/2\}\cap U_1.\]
Then $W_t$ is an open cover of the compact set $\supp(t\mapsto F(t,\cdot))$, so there exists a finite subcover $W_{t_1},\ldots W_{t_N}$.  Let $\{\psi_i\}_{i=1}^N$ be a partition of unity subordinate to this cover.
Since $\sum\psi_i(t)\leq 1$ for all $t\in \T$, 
\begin{align*}
\|\sum_{i=1}^N\psi_i(t)F(t_i,\cdot)-F(t,\cdot)\|_\infty &=\|\sum_{i=1}^N \psi_i(t) F(t_i,\cdot)-\sum_{i=1}^N \psi_i(t)F(t,\cdot)\|_{\infty}\\
&\leq \sum_{i=1}^N\psi_i(t)\|F(t_i,\cdot)-F(t,\cdot)\|_{\infty}<\frac{\epsilon}{2}.
\end{align*}

 For $\gamma\in \bigcup_{i=1}^N \supp(F(t_i,\cdot))\subset U_2$, let $U_\gamma$ be the open set 
\[
U_\gamma:=\{\eta\in G:
|F(t_i,\gamma)-F(t_i,\eta)|<\epsilon/2\text{\ for  $1\leq i\leq N$}\}\cap U_2.
\]
Since $\bigcup_{i=1}^N \supp(F(t_i,\cdot))$ is compact there exists a finite subcover $U_{\gamma_1},\ldots, U_{\gamma_M}$. Let  $\{f_j\}_{j=1}^M$ be a partition of unity subordinate to this subcover.   For $\gamma\in G$ and each $i\in\{1,\ldots, N\}$ we have
\begin{align*}
\big|\sum_{j=1}^M f_j(\gamma)F(t_i,\gamma_j)-F(t_i,\gamma)\big| &=\big|\sum_{j=1}^M f_j(\gamma)F(t_i,\gamma_j)-\sum_{j=1}^M f_j(\gamma)F(t_i,\gamma)\big| \\
&\leq \sum_{j=1}^M f_j(\gamma)\big|F(t_i,\gamma_j)-F(t_i,\gamma)\big|<\frac{\epsilon}{2}.
\end{align*}
Now set $F_\epsilon:=\sum_{i,j=1}^{M,N} F(t_i,\gamma_j)\psi_i\otimes f_j$, and  note that $\supp F_\epsilon$ is contained in $U_1\times U_2$ by construction. We have
\begin{align*}
&\|F_\epsilon-F\|_\infty=\sup_{(t,\gamma)}\{\big|\sum_{i,j}F(t_i,\gamma_j)\psi_i(t)f_j(\gamma)-F(t,\gamma)\big|\}\\
&\leq \sup_{(t,\gamma)}\big\{\big|\sum_{i,j}\psi_i(t)f_j(\gamma)F(t_i,\gamma_j)-\sum_{i}\psi_i(t)F(t_i,\gamma)+\sum_{i}\psi_i(t)F(t_i,\gamma)-F(t,\gamma)\big|\big\}\\
&<\sup_{(t,\gamma)}\big\{\sum_{i}\psi_i(t)\big|\sum_{j}f_j(\gamma)F(t_i,\gamma_j)-F(t_i,\gamma)\big|+\frac{\epsilon}{2}\big\}<\frac{\epsilon}{2}+\frac{\epsilon}{2}=\epsilon.\end{align*}
We have now shown that $\spn\{\psi\otimes f:\psi\in C(\T), f\in C_c(G)\}$ is dense in $C_c(G^\omega)$ in the inductive limit topology.  Thus it follows from the Stone-Weierstrass theorem that  $\spn\{s^m\otimes f:m\in \Z, f\in C_c(G)\}$ is dense in $C_c(G^\omega)$.\end{proof}

Lemmas~\ref{lem-Im} and \ref{lem: density of poly} give:

\begin{prop}\label{prop dec E} Let $G^\omega$ be the groupoid extension associated to $(G,\omega)$. Then $C^*(G^\omega)=\oplus_{n\in\Z}I_n$.
\end{prop}

\begin{proof}[Proof of Theorem~\ref{thm-dec}]
Both $\chi_{n}:I_{n}\to C^*(G^\omega, n)$ and $\phi_{n}:C^*(G^\omega,n)\to C^*(G, \omega^n)$ are isomorphisms by Lemmas~\ref{lem:chi iso} and  \ref{lem-algberas}\eqref{item:isotoRenault},  so $\Upsilon_n|_{I_{n}}=\phi_n\circ \chi_n|_{I_n}$ is an isomorphism of $I_{n}$ onto $C^*(G, \omega^n)$. But by Lemma~\ref{lem:chi iso}\eqref{item3:chi iso}, $\chi_{n}(I_m)=\{0\}$ if $n\neq m$, so Theorem~\ref{thm-dec} follows from Proposition~\ref{prop dec E}.  
\end{proof}

We now show that Theorem~\ref{thm-dec} leads to a general framework for deducing results about  twisted groupoid $C^*$-algebras from untwisted ones.  The basic idea is that many properties of principal groupoids  are shared with their extensions by $\T$. 
We start with a general lemma. The \emph{stabilizer subgroupoid}  of a groupoid $G$ is $\{\gamma\in G:r(\gamma)=s(\gamma)\}$ and $A_u:=\{\gamma\in G: r(\gamma)=u=s(\gamma)\}$ is  the \emph{stability subgroup} at $u$. 

\begin{lem}\label{lem-new}
Let  $G$ be a groupoid and  $G^\omega$ be the extension associated to $(G,\omega)$. Let $A$ and $A^\omega$ be the respective stabilizer subgroupoids of $G$ and $G^\omega$.  Then the map $(t,\gamma)\mapsto [\gamma]$ induces a homeomorphism and  isomorphism of  $G^\omega/A^\omega$ onto $G/A$. 
\end{lem}

\begin{proof} The stability subgroups of $G^\omega$ are $\T\times_\omega A_u$ where $A_u$ is the stability subgroup of $G$ at $u$.  Thus $A^\omega=\cup_{u\in G^{(0)}}\T\times_\omega A_u=\T\times_\omega A$.

  We will show that the map $f:(t,\gamma)\mapsto [\gamma]$ induces a homeomorphism and  isomorphism of $G^\omega/A^\omega$ onto $G/A$. Certainly $f$ is a groupoid morphism, and is continuous and surjective.  If $f(t,\gamma)=f(s,\delta)$ then there exists $\alpha\in A_{s(\gamma)}$ such that $\gamma=\delta\alpha$.  Then $(t,\gamma)=(s,\delta)(s^{-1}t\overline{\omega(\delta, \alpha)}, \alpha)$, and $(s^{-1}t\overline{\omega(\delta, \alpha)}, \alpha)\in A^\omega$. Hence $[(t,\gamma)]=[(s,\delta)]$.  So $f$ induces a continuous bijection $\tilde f:G^\omega/A^\omega\to  G/A$. Similarly, the function $g:G\to G^\omega/A^\omega$ defined by $g(\gamma)=[(1,\gamma)]$ induces a continuous bijection $\tilde g:G/A\to G^\omega/A^\omega$, and it is easy to check that $\tilde g$ is the inverse of $\tilde f$.  Thus $\tilde f$ is a homeomorphism.  
\end{proof}

\begin{prop}\label{prop-similar} Let $G$ be a principal groupoid and  let $G^\omega$ be the extension associated to a continuous $2$-cocycle $\omega: G^{(2)}\to\T$. Then 
\item \begin{enumerate} 
\item\label{prop-similar1} $C^*(G)$ has continuous trace if and only if $C^*(G^\omega)$ has continuous trace; and 
\item\label{prop-similar2} if $C^*(G)$ has continuous trace then so does $C^*(G,\omega)$.
\end{enumerate} 
\end{prop}

\begin{proof} \eqref{prop-similar1} First suppose that $C^*(G)$ has continuous trace. Since $G$ is principal, \cite[Theorem~2.3]{MW90} implies that $G$ is a proper groupoid. Now consider $G^\omega$: since $G$ is principal the stability subgroups of $G^\omega$ are  $\T\times \{u\}$ where $u\in G^{(0)}$, and the stabilizer subgroupoid is $A^\omega=\T\times G^{(0)}$. In particular,  the stability subgroups of $G^\omega$ are all abelian and $u\mapsto \T\times \{u\}$ is continuous in the Fell topology on the set of closed subgroups of $G$.  By Lemma~~\ref{lem-new}, the quotient groupoid $G^\omega/A^\omega$ and $G$ are homeomorphic, and hence $G^\omega/A^\omega$ is proper.  Now $C^*(G^\omega)$ has continuous trace by \cite[Theorem~1.1]{MRW96}.

Conversely, suppose  $C^*(G^\omega)$ has continuous trace. By Theorem~\ref{thm-dec}, $C^*(G)=C^*(G,\omega^0)$ is a direct summand of $C^*(G^\omega)$. Hence $C^*(G)$ has continuous trace by \cite[Proposition~6.2.10]{Pedersen}. 

\eqref{prop-similar2}  Suppose that $C^*(G)$ has continuous trace. Then $C^*(G^\omega)$ has continuous trace by \eqref{prop-similar1}. By Theorem~\ref{thm-dec}, $C^*(G,\omega)$ is a direct summand of $C^*(G^\omega)$,   and hence $C^*(G,\omega)$ has continuous trace as well.
\end{proof}

Many properties are shared by $G$ and $G^\omega$: having a Haar system,  being Cartan,  proper or integrable, and any topological property of the orbit spaces. This gives the proposition below. Example~\ref{ex-stuffed} below shows that we cannot expect to extend Propositions~\ref{prop-similar} and \ref{prop-similaragain} to non-principal groupoids $G$.

\begin{prop}\label{prop-similaragain} Let $G$ be a principal groupoid and   let $G^\omega$ be the extension associated to a continuous $2$-cocycle $\omega: G^{(2)}\to\T$.
\begin{enumerate}
\item\label{prop-similaragain1} $C^*(G)$ is a Fell algebra if and only if  $C^*(G^\omega)$ is a Fell algebra. 
If $C^*(G)$ is a Fell algebra then  so is   $C^*(G,\omega)$.
\item\label{prop-similaragain2} $C^*(G)$ has bounded trace  if and only if $C^*(G^\omega)$   has bounded trace.
If $C^*(G)$  has bounded trace  then  so does  $C^*(G,\omega)$.
\item\label{prop-similaragain3} $C^*(G)$ is liminal if and only if $C^*(G^\omega)$ is liminal.
If $C^*(G)$ is liminal then so is $C^*(G,\omega)$;
\item\label{prop-similaragain4} $C^*(G)$ is postliminal if and only if $C^*(G^\omega)$ is postliminal.
If $C^*(G)$ is postliminal then so is $C^*(G,\omega)$.
\end{enumerate}
\end{prop}
\begin{proof} Since $G$ is principal the stability subgroups of $G^\omega$ are $\T\times\{u\}$ where $u\in G^{(0)}$; in particular they are abelian and vary continuously.

\eqref{prop-similaragain1} and \eqref{prop-similaragain2} We can proceed as in the proof of Proposition~\ref{prop-similar} replacing \cite[Theorem~1.1]{MRW96} with \cite[Theorem~6.5]{CH10} and  \cite[Theorem~6.4]{CH10}, respectively.

\eqref{prop-similaragain3}
First suppose that $C^*(G)$ is liminal. Since the stability subgroups of $G$ are trivial, the orbit space $G^{(0)}/G$ is $T{}_{1}$ by \cite[Theorem~6.1]{Cla07}. But the orbit space of $G^\omega$ is homeomorphic to the orbit space of $G$ via $[(1,u)]\mapsto [u]$, hence is  $T{}_{1}$ as well. Since the stability subgroups of $G^\omega$ are amenable and liminal, $C^*(G^\omega)$ is liminal by  \cite[Theorem~6.1]{Cla07}.  

Second, suppose  that $C^*(G^\omega)$ is liminal. By Theorem~\ref{thm-dec}, $C^*(G)=C^*(G,\omega^0)$ is a direct summand of $C^*(G^\omega)$. Hence $C^*(G)$ is liminal by \cite[Proposition~6.2.9]{Pedersen}. This gives the first statement of  \eqref{prop-similaragain3}. 

Finally,  suppose that $C^*(G)$ is liminal. Then $C^*(G^\omega)$ is liminal.  By Theorem~\ref{thm-dec}, $C^*(G,\omega)$ is a direct summand of $C^*(G^\omega)$, and hence    $C^*(G,\omega)$ must be liminal as well. This gives the second statement of  \eqref{prop-similaragain3}. 

\eqref{prop-similaragain4} Theorem~7.1 of \cite{Cla07} says that the groupoid $C^*$-algebra of a groupoid with amenable stability subgroups is postliminal if and only if the orbit space is  T${}_0$
and the stability subgroups are postliminal.  So  \eqref{prop-similaragain4}  follows as above using  \cite[Theorem~7.1]{Cla07} in place of  \cite[Theorem~6.1]{Cla07}.
\end{proof}

\begin{ex}\label{ex-stuffed} When the groupoid $G$ is not principal, the stability subgroups of $G^\omega$ can easily fail to be abelian, liminal or postliminal even if the stability subgroups of $G$ are abelian, liminal or postliminal, respectively. Thus the theorems used to prove Propositions~\ref{prop-similar} and \ref{prop-similaragain}, such as \cite[Theorem~1.1]{MRW96} and \cite[Theorem~6.1]{Cla07}, do not apply.   
The following is an example of a group $G$ and a $2$-cocycle $\omega$ such that  $C^*(G)$ has continuous trace but $C^*(G,\omega)$ and $C^*(G^\omega)$ are not even postliminal. Thus we cannot expect an analog of Proposition~\ref{prop-similar} when the groupoid $G$ is not principal.

Let $\theta\in(0,1)$ be irrational and define $\omega:\Z^2\times\Z^2\to\T$ by $\omega((m_1,m_2),(n_1, n_2))=e^{-2\pi im_1n_2\theta}$.  The twisted group $C^*$-algebra $C^*(\Z^2,\omega)$ is isomorphic to the irrational rotation algebra $A_\theta=C(\T)\rtimes\Z$ (see, for example, \cite[pp. 21-22]{DR}). Since $\theta$ is irrational the orbit space $\T/\Z$ is not T${}_0$, and hence $A_\theta$ is not postliminal by \cite[Theorem~3.3]{Gootman}. Thus $C^*(\Z^2,\omega)$ is not postliminal. By Theorem~\ref{thm-dec},  $C^*(\Z^2,\omega)$ is a summand of $C^*((\Z^2)^{\omega})$, so $C^*((\Z^2)^{\omega})$ is not postliminal either. 
Thus $C^*(\Z^2,\omega)$ and $C^*((\Z^2)^{\omega})$ are not postliminal even though $C^*(\Z^2)\cong C(\T^2)$ has continuous trace. 
\end{ex}


\section{A reduced version of the decomposition theorem} The goal of this section is to prove a version of Theorem~\ref{thm-dec} for reduced crossed products.  Let $E$ be a second-countable locally compact Hausdorff groupoid with a Haar system $\beta$,  and $\sigma:E^{(2)}\to\T$ be a continuous $2$-cocycle.   Let  $u\in E^{(0)}$. The \emph{left-regular representation} $\Pi^u$ is the representation of $C_c(E,\sigma)$ on $L^2(E,\beta_u)$ characterized by
\begin{equation}\label{left reg}
\rinner{\Pi^u(f)\xi}{\zeta}{}=\int_E\int_E f(\gamma\eta)\xi(\eta\inv)\overline{\zeta(\gamma)}\sigma(\gamma\eta,\eta\inv)\, d\beta^u(\gamma)\, d\beta_u(\eta)
\end{equation}
for $f\in C_c(E,\sigma)$ and $\xi,\zeta\in L^2(E,\beta_u)$.  Since $\Pi^u$ is continuous in the inductive limit topology, it extends to a representation of $C^*(E,\sigma)$. The reduced $C^*$-algebra $C^*_r(E,\sigma)$  of $(E,\sigma)$ is the completion of $C_c(E,\sigma)$ with respect to  the norm $\|f\|_r=\sup_{u\in E^{(0)}}\{\|\Pi^u(f)\|\}$.  Alternatively, $C^*_r(E,\sigma)=C^*(E,\sigma)/I$ where $I=\bigcap_{u\in E^{(0)}}\ker(\Pi^u)$; we write $q=q_{_{E}}$ for the quotient map.

\begin{thm}\label{thm-red dec} Let $G$ be a second-countable, locally compact Hausdorff groupoid with a Haar system $\lambda$.  Let $\omega:G^{(2)}\to\T$ be a continuous 2-cocycle and ${G^\omega}$ the extension associated to $(G,\omega)$.  Let $\Upsilon:C^*({G^\omega})\to\oplus_{n\in\Z}C^*(G, \omega^n)$ and $\Upsilon_n:C^*({G^\omega})\to      C^*(G, \omega^n)$ be as in Theorem~\ref{thm-dec}. Then there exists a homomorphism $\Omega_n: C_r^*({G^\omega})\to C^*_r(G,\omega^n)$ such that the following diagram 
\begin{equation}\label{red diagram}\xymatrix{C^*({G^\omega}) \ar[r]^{\Upsilon_n} \ar[d]^{q_{_{G^\omega}} }& C^*(G,\omega^n) \ar[d]^{q_{_{G,n}}}\\ C_r^*({G^\omega})\ar[r]^{\Omega_n} &C^*_r(G,\omega^n)}
\end{equation}  
commutes.         
Furthermore,  the map  $\Omega:C^*_r({G^\omega})\to \oplus_{n\in\Z} C^*_r(G, \omega^n)$, defined by $a\mapsto (\Omega_n(a))$, is an isomorphism.
\end{thm}

For $n\in \Z$  and $u\in \units$,   we write $L^u_n$ for the left-regular representation  of $C^*(G,\omega^n)$  on $L^2(G,\lambda_u)$  
and $R^u$ for the left-regular representation of $C^*({G^\omega})$ on $L^2({G^\omega}, \tau\times\lambda_u)$; both are characterized by \eqref{left reg}.

\begin{lem}\label{lem: lots of stuff} Let $G^\omega$ be the extension associated to $(G,\omega)$. Let $u\in \units$. For $n\in \Z$ define 
$$\mathcal{H}_n^u:=\overline{\spn}\{s^{-n}\otimes \xi:\xi\in C_c(G)\}\subset L^2({G^\omega},\tau\times\lambda_u).$$

\begin{enumerate}
\item\label{stuff a} If $m\neq n$ then $\mathcal{H}^u_m$ is orthogonal to $\mathcal{H}^u_n$.
\item\label{stuff b} There is a unitary $V_n:L^2(G,\lambda_u)\to \mathcal{H}_n^u$ such that 
\begin{equation}\label{Vnformula}
V_n(\xi)=s^{-n}\otimes\xi\quad\text{for}\quad\xi\in C_c(G).
\end{equation}
\item\label{stuff c}There is a unitary $V:\oplus_{n\in\Z}L^2(G,\lambda_u)\to L^2({G^\omega},\tau\times\lambda_u)$ characterized by  $$V((\xi_n))=\oplus_{n\in \Z} V_n(\xi_n)\quad\text{ for }\quad \xi_n\in C_c(G).$$

\item\label{stuff d}For $n\in \Z$ let $L_n^u:C^*(G,\omega^n)\to B(L^2(G,\lambda_u))$  and $R^u:C^*({G^\omega})\to B(L^2({G^\omega},\tau\times\lambda_u))$ be the respective left-regular representations, and set $L^u=\oplus_{n\in\Z}L_n^u$. Then
$$V(L^u\circ\Upsilon(a))V^*=R^u(a)\quad \text{for all}\quad a\in C^*({G^\omega}).$$
\end{enumerate}
\end{lem}

\begin{proof}

We compute:
\begin{align}\rinner{r^{-m}\otimes\xi}{r^{-n}\otimes\zeta}{L^2({G^\omega})}&=\int_G\int_\T r^{-m}\otimes\xi(t,\gamma)\overline{r^{-n}\otimes\zeta(t,\gamma)}\, dt \, d\lambda_u(\gamma)\nonumber\\
&=\int_G\int_\T t^{n-m}\xi(\gamma)\overline{\zeta(\gamma)}\, dt\, d\lambda_u(\gamma)\nonumber\\
&=\rinner{\xi}{\zeta}{L^2(G)}\delta_{m,n}.\label{eq:  rel inner}\end{align}
Now~\eqref{eq:  rel inner} implies, first, that $\mathcal{H}^u_m$ is orthogonal to $\mathcal{H}^u_n$, and second, that there is an isometry $V_n$ satisfying  \eqref{Vnformula}. By definition of $\mathcal{H}^u_n$, $V_n$ is onto and  hence is unitary.  This gives \eqref{stuff a} and \eqref{stuff b}.

By Lemma~\ref{lem: density of poly},  $\spn\{r^m\otimes \xi:m\in \Z, \xi\in C_c(G)\}$ is dense in $C_c({G^\omega})$ in the inductive limit topology, and hence it is dense in $L^2({G^\omega},\tau\times \lambda_u)$ as well.  Now~\eqref{stuff c} follows from \eqref{stuff a} and \eqref{stuff b}.

For \eqref{stuff d}, let $m, n\in \Z$, $\xi,\zeta\in C_c(G)$ and $F\in C_c({G^\omega})$. Then, using Fubini's Theorem several times,  
\begin{align*}
&\rinner{R^u(F)(r^{-m}\otimes\xi)}{r^{-n}\otimes\zeta}{L^2(G^\omega)}\nonumber\\
&=\int_G\int_\T\int_G\int_\T F\big((t,\gamma)(s,\eta)  \big) (r^{-m}\otimes\xi)\big( (s,\eta)^{-1} \big) \overline{r^{-n}\otimes\zeta(t,\gamma)}\, ds \, d\lambda^u(\eta) \, dt \, d\lambda_u(\gamma)\\
&=\int_G\int_G\int_\T\int_\T F(st\omega(\gamma,\eta),\gamma\eta)s^{m}\omega(\eta,\eta\inv)^m\xi(\eta\inv)t^{n}\overline{\zeta(\gamma)}\, ds \, dt \, d\lambda^u(\eta)\, d\lambda_u(\gamma)\nonumber\\
\intertext{and, replacing $s$ with $st\inv\overline{\omega(\gamma,\eta)}$, gives}
&=\int_G\int_G\int_\T\int_\T F(s,\gamma\eta)s^{m}t^{-m}\overline{\omega(\gamma,\eta)}^{m}\omega(\eta,\eta\inv)^m\xi(\eta\inv)t^{n}\overline{\zeta(\gamma)}\, ds \, dt \, d\lambda^u(\eta)\, d\lambda_u(\gamma)\nonumber\\
\intertext{which, because  $\omega(\gamma,\eta)\omega(\gamma\eta,\eta\inv)=\omega(\eta,\eta\inv)\omega(\gamma,\eta\eta\inv)=\omega(\eta,\eta\inv)$, becomes}
&=\int_G\int_G\bigg(\int_\T F(s,\gamma\eta)s^{m}\, ds \bigg)\biggl(\int_\T t^{n-m} \, dt\biggr) \omega(\gamma\eta,\eta\inv)^{m}\xi(\eta\inv)\overline{\zeta(\gamma)} \, d\lambda^u(\eta)\, d\lambda_u(\gamma)\nonumber\\
&=\delta_{m,n}\int_G\int_G \Upsilon_{m}(F)(\gamma\eta)\xi(\eta\inv)\overline{\zeta(\gamma)}\omega(\gamma\eta,\eta\inv)^{m}\, d\lambda^u(\eta)\, d\lambda_u(\gamma)\nonumber\\
&=\delta_{m,n}\rinner{L_{m}^u(\Upsilon_{m}(F))\xi}{\zeta}{L^2(G)}.\label{eq: rs1l}
\end{align*}
Since $C_c(G^\omega)$ is dense in $C^*(G^\omega)$,  it follows that for $a\in C^*(G^\omega)$
\begin{align*}\rinner{R^u(a)\bigg(\sum r^{-m}\otimes \xi_m\bigg)}{\sum r^{-n}\otimes\zeta_n}{L^2(G^\omega)}&=\sum_{m,n}\rinner{L_{m}^u(\Upsilon_{m}(a))\xi_m}{\zeta_n}{L^2(G)}\delta_{m,n}\\
&=\sum_n\rinner{L_{n}^u(\Upsilon_{n}(a))\xi_n}{\zeta_n}{L^2(G)}
.\end{align*}

  So for $x=\sum s^m\otimes \xi_m$, $y=\sum s^n\otimes\zeta_n$ we have
\begin{align*}
\rinner{R^u(a)x}{y}{L^2(G^\omega)}
&= \sum_n\rinner{L^u_n(\Upsilon_{n}(a))\xi_n}{\zeta_n}{L^2(G)}
=\rinner{L^u(\Upsilon(a))V^*x}{V^*y}{\oplus_{n\in \Z} L^2(G)},
\end{align*}
and then~\eqref{stuff d} follows because the set of such $x,y$ is dense in $L^2({G^\omega},\tau\times\lambda_u)$.
\end{proof}

\begin{proof}[Proof of Theorem~\ref{thm-red dec}] By  Lemma \ref{lem: lots of stuff}\eqref{stuff d}, we have $\ker(R^u)\subset\ker(L^u_n\circ \Upsilon_{n})$ for all $n$.  Since this holds for all $u\in \units$, $\ker(q_{_{G^\omega}})\subset\ker(q_{_{G,n}})$. Thus the map $q_{_{G,n}}\circ \Upsilon_n$ induces a homomorphism $\Omega_n$ such that the diagram \eqref{red diagram}
commutes. 

To see that $\Omega=(\Omega_n)$ is isometric, recall from Proposition~\ref{prop dec E} that $C^*(G^\omega)=\oplus_{m\in\Z}I_m$ and  let $a=(a_n)\in C^*(G^\omega)$ where $a_n\in I_n$. Using first Lemma~\ref{lem: lots of stuff}\eqref{stuff d}, and second,   $\Upsilon_n=\Upsilon|_{I_n}$ and $L^u=\oplus_n L^u_n$, we get 
$$\|R^u(a)\|=\|L^u(\Upsilon(a))\|=\max_n\|L^u_n(\Upsilon_n(a_n))\|.$$
Since this holds for all $u\in\units$,
\begin{align*}\|q_{_{G^\omega}}(a)\|_{C^*_r(G^\omega)}&=\max_n\|q_{_{G,n}}(\Upsilon_n(a_n))\|_{C^*_r(G,\omega^n)}=\max_n\|\Omega_n(q_{_{G^\omega}}(a_n))\|_{C^*_r(G,\omega^n)}\\
&=\|\Omega(q_{_{G^\omega}}(a))\|_{C^*_r(G,\omega^n)}.
\end{align*}
Hence $\Omega$ is isometric.  That $\Omega$ is surjective follows from the commutativity of the diagram since $\Upsilon=(\Upsilon_n)$ and the quotient maps are surjective. Thus $\Omega$ is an isomorphism.
\end{proof}

\begin{cor}\label{cor amen}
Let $G^\omega$ be the extension associated to $(G,\omega)$. If $G$ is amenable, then $C^*(G,\omega)=C^*_r(G,\omega)$.
\end{cor}

\begin{proof}
Let $j:G^\omega\to G$ be the quotient map.  Then $\ker j=\T\times \units$ is amenable.  Since $G$ is  amenable, 
Proposition~5.1.2 of \cite{AR00} implies that $G^\omega$ is amenable.  By Theorems~\ref{thm-dec} and \ref{thm-red dec} we have
\[
\oplus_{n\in\Z} C^*_r(G, \omega^n)\cong C^*_r(G^\omega)= C^*(G^\omega)\cong \oplus_{n\in\Z} C^*(G, \omega^n).
\]
By the commutativity of \eqref{red diagram} the summands corresponding to $n=1$ match up,  so the result follows.
\end{proof}


\section{Actions of proper groupoids and fixed-point algebras}\label{sec-app}

Let $G$ be a principal proper groupoid.  Then $G^\omega=\T\times_\omega G$ is also proper.  There is an action $\lt$ of $G^\omega$ on $C_0(\units)$ defined  by 
\[\lt_\gamma(f)(v)=f(s_{G^\omega}(\gamma))\quad\text{ for $f\in C_0(\units)$  and $\gamma\in G^\omega$ with $r_{G^\omega}(\gamma)=v$.}
\]
 Since $G^\omega\backslash \units = G\backslash \units$, \cite[Proposition 2.2]{MW90} implies that $C_0(G^\omega\backslash \units)$ is Morita equivalent to $C^*(G)$.  Theorem 3.9 of \cite{mep09} implies that $C_0(G^\omega\backslash \units)$ is Morita equivalent to an ideal $I$ of $C^*(G^\omega)$.  In the following proposition we reconcile these two results by using the decomposition of $C^*(G^\omega)$ into the direct sum $\oplus_{n\in\Z} C^*(G,\omega^n)$ to identify the ideal $I$ with the summand corresponding to $n=0$. 

\begin{prop}\label{cor-notsat} 
Let $G$ be a principal and proper groupoid.  Then the generalized fixed-point algebra $C_0(\units)^{\lt}=C_0({G^\omega}\backslash \units)$ is Morita equivalent  to the direct summand $C^*(G)=C^*(G,\omega^{0})$ of $C^*(G^\omega)$. \end{prop}

\begin{proof}
Theorem~3.9 of \cite{mep09} says that there is a $C^*$-subalgebra $I$ of the reduced groupoid crossed product $C_0(\units)\rtimes_{\lt, r}{G^\omega}$ that is Morita equivalent to a generalized fixed-point algebra $C_0(\units)^{\lt}$ of $(C_0(\units), G^\omega,\lt)$. In our special case where the groupoid acts properly on its unit space, $I$ is an ideal  by Remark~4.14 of \cite{mep09}.  By Proposition~4.1 of \cite{mep09}, $C_0(\units)^{\lt}=C_0({G^\omega}\backslash \units)$.  Combining \cite[Corollary~2.1.7 and Proposition~3.3.5]{AR00} gives that ${G^\omega}$ is measurewise amenable, and hence $C_0(\units)\rtimes_{\lt, r}{G^\omega}=C_0(\units)\rtimes_{\lt}{G^\omega}$ by \cite[Proposition~6.1.10]{AR00}.  By \cite[Remark~4.22]{Goehle09}, $C_0(\units)\rtimes_{\lt}{G^\omega}$ is isomorphic to $C^*({G^\omega})$.   So $I$ is an ideal in $C^*({G^\omega})$ that is Morita equivalent to $C_0({G^\omega}\backslash \units)$; it remains to identify the ideal $I$.  

The imprimitivity bimodule implementing the Morita equivalence is a completion of $C_c(\units)$ with respect to the left inner  product given by 
\begin{equation}\label{eq-compare}
\linner{I}{f}{g}(t,\gamma) =f(r_{G^\omega}(t,\gamma))\overline{g(s_{G^\omega}(t,\gamma))}=f(r_{G}(\gamma))\overline{g(s_G(\gamma))}
\end{equation}
for $f,g\in C_c(\units)$ and $(t,\gamma)\in \T\times_\omega G$.  The point is that the inner product is independent of  $t$. Thus $I$ is an ideal of $C^*({G^\omega},0)\cong C^*(G)$. When we apply Theorem~3.9 of \cite{mep09} to the action $\lt$ of $G$ on $C_0(G^{(0)})$ we obtain a Morita equivalence based on $C_c(G^{(0)})$ between an ideal $J$ of $C^*_r(G)$ and $C_0(G\backslash G^{(0)})$, with left inner product given by \begin{equation}\label{eq-compare2}\linner{J}{f'}{g'}(\gamma) =f'(r_G(\gamma))\overline{g'(s_G(\gamma))}\end{equation} 
for $f', g'\in C_c(G^{(0)})$ and $\gamma\in G$. Note that $C^*(G)=C^*_r(G)$ by amenability.
Comparing \eqref{eq-compare} and \eqref{eq-compare2} shows $I=J$.  Finally, by \cite[Theorem~5.9]{mep09}, $J=C^*(G)$.
\end{proof}

The action $\lt$ is called \emph{saturated} if the ideal $I$ of $C_0(\units)\rtimes_{\lt, r}{G^\omega}$ is in fact $C_0(\units)\rtimes_{\lt, r}{G^\omega}$. We note that in the situation of Proposition~\ref{cor-notsat} the action is very far away from being saturated since it is just one summand in $\oplus_{n\in\Z}C^*(G, \omega^n)$.




\providecommand{\bysame}{\leavevmode\hbox to3em{\hrulefill}\thinspace}
\providecommand{\MR}{\relax\ifhmode\unskip\space\fi MR }
\providecommand{\MRhref}[2]{%
  \href{http://www.ams.org/mathscinet-getitem?mr=#1}{#2}
}
\providecommand{\href}[2]{#2}

\end{document}